\documentclass{article}
\usepackage{amsmath,amsfonts,amsthm,amssymb}

\DeclareMathSymbol\nullset{\mathord}{AMSb}{"3F}

\begin{document}
\newtheorem{theorem}{Theorem}[section]
\newtheorem{corollary}{Corollary}[section]
\newtheorem{lemma}{Lemma}[section]
\newtheorem{proposition}{Proposition}[section]
{\newtheorem{definition}{Definition}[section]}

\title{Examples of associative algebras for which the $T$-space of central polynomials is not finitely based}
\author{C. Bekh-Ochir and S. A. Rankin
  }

\maketitle
\newcounter{parts}
\def\set#1\endset{\{\,#1\,\}}
\def\gen#1{\left<{#1}\right>}
\def\rest#1{{}_{{}_{#1}}}
\def\com#1,#2{[{#1},{#2}]}
\def\choice#1,#2{\binom{#1}{#2}}
\def\kx{k\langle X\rangle}
\def\kzerox{k_0\langle X\rangle}
\def\konex{k_1\langle X\rangle}
\def\unitgrass{G}
\def\nonunitgrass{G_0}
\def\finiteunitgrass#1{G(#1)}
\def\finitenonunitgrass#1{G_0(#1)}
\def\siderovset{SS}
\def\boundedsiderovset{BSS}
\def\specialsiderovset{SPSS}
\def\finitesiderovsetone#1{SS'(#1)}
\def\finitesiderovsettwo#1{SS''(#1)}
\def\boundedsiderovset{BSS}
\def\bss#1{BSS(#1)}
\def\lend(#1){lend(#1)}
\def\lbeg#1{lbeg(#1)}
\def\endof(#1){end(#1)}
\def\beg(#1){beg(#1)}
\def\finitesidset#1{SS(#1)}
\let\cong=\equiv
\def\mod#1{\,\,(\text{mod}\,#1)}
\def\mz(#1){M_0(#1)}
\def\wt(#1){\text{wt}(#1)}
\def\dom(#1){\text{dom}(#1)}
\def\from{\mkern2mu{:}\mkern2mu}
\def\zplus{\mathbb{Z}^+}

\begin{abstract}
In 1988 (see \cite{Ok}), S. V.  Okhitin proved that for any field $k$ of characteristic zero, the $T$-space 
$CP(M_2(k))$ is finitely based, and he raised the question as to whether $CP(A)$ is finitely based for every 
(unitary) associative algebra $A$ for which $0\ne T(A)\subsetneq CP(A)$.
V. V. Shchigolev (see \cite{Sh}, 2001) showed that for any field of characteristic zero, every $T$-space of $\kzerox$
is finitely based, and it follows from this that every $T$-space of $\konex$ is also finitely based. This more than
answers Okhitin's question (in the affirmative) for fields of characteristic zero.

For a field of characteristic 2, the infinite-dimensional Grassmann algebras, unitary and nonunitary, are commutative
and thus the $T$-space of central polynomials of each is finitely based. 

We shall show in the following that if $p>2$ and $k$ is an arbitrary field of characteristic $p$, then neither
$CP(\nonunitgrass)$ nor $CP(\unitgrass)$ is finitely based, thus providing a negative answer to
Okhitin's question.
\end{abstract}



\section{Introduction and preliminaries}
 Let $k$ be a field of characteristic $p$ and let $X$ be a countably infinite set, 
 say $X=\set x_i\mid i\ge 0\endset$. Then $\kzerox$ denotes the free associative 
 $k$-algebra on $X$, while $\konex$ denotes the free unitary associative
 $k$-algebra on $X$. 

 Let $A$ denote any associative $k$-algebra. 
 For any $H\subseteq A$, $\gen{H}$ shall denote the linear subspace of $A$ spanned by $H$.
 Any linear subspace of $A$ that is invariant under every endomorphism of $A$ is
 called a $T$-space of $A$, and if a $T$-space happens to also be an ideal of $A$, then it
 is called a $T$-ideal of $A$. For $H\subseteq A$, the smallest $T$-space of $A$ that contains $H$
 shall be denoted by $H^S$, while the smallest $T$-ideal of $A$ that contains $H$ shall be
 denoted by $H^T$. If $V\subseteq A$ is a $T$-space and there exists finite $H\subseteq A$ such that
 $V=H^S$, then $V$ is said to be finitely based. In this article, we shall deal only with $T$-spaces and
 $T$-ideals of $\kzerox$ and $\konex$. Occasionally, we shall consider $H\subseteq \kzerox\subseteq \konex$,
 and we may wish to have notation for both the $T$-space generated by $H$ in $\kzerox$ and the $T$-space
 generated by $H$ in $\konex$ so that both could appear in the same discussion. Accordingly, we shall
 write $H^{S_0}$ to denote the $T$-space of $\kzerox$ that is generated by $H$, and let $H^S$ denote the
 $T$-space of $\konex$ that is generated by $H$. Similarly, we may use $H^{T_0}$ to denote the $T$-ideal of
 $\kzerox$ that is generated by $H$, and $H^T$ for the $T$-ideal of $\konex$ that is generated by
 $H$.

 A nonzero element $f\in \kzerox$ is called an {\em identity} of $A$
 if $f$ is in the kernel of every homomorphism from $\kzerox$ to $A$ (every unitary homomorphism from $\konex$ if $A$ is unitary). 
 The set of all identities of $A$, together with 0, forms a $T$-ideal of $\kzerox$ (and of $\konex$ if $A$ is unitary), 
 denoted by $T(A)$. An element $f\in \kzerox$ is called a {\it central polynomial} of $A$ 
 if $f\notin T(A)$ and the image of $f$ under any homomorphism from $\kzerox$ (unitary homomorphism from $\konex$ if $A$ is unitary)
 belongs to $C_{A}$, the centre of $A$. The $T$-space of $\kzerox$ (or of $\konex$ if $A$ is unitary)
 that is generated by the set of all central polynomials of $A$ is denoted by $CP(A)$.

 Let $\unitgrass$ denote the (countably) infinite dimensional unitary
 Grassmann algebra over $k$, so there exist $e_i\in \nonunitgrass$, $i\ge 1$, 
 such that  for all $i$ and $j$, $e_ie_j=-e_je_i$, $e_i^2=0$, and
 $\mathcal{B}=\set e_{i_1}e_{i_2}\cdots e_{i_n}\mid n\ge 1,
 i_1<i_2<\cdots i_n\endset$, together with $1$, forms a linear basis for
 $G$. Let $E$ denote the set $\set e_i\mid i\ge1\endset$. The subalgebra of $\unitgrass$ with linear basis $\mathcal{B}$ is
 the infinite dimensional nonunitary Grassmann algebra over $k$, and is
 denoted by $\nonunitgrass$. 

 It is well known that $T^{(3)}$, the $T$-ideal of $\konex$ generated  by $\com
 {\com x_1,{x_2}},{x_3}$, is contained in $T(\unitgrass)$. For convenience, 
 we shall write $\com x_1,{x_2,x_3}$ for  $\com {\com x_1,{x_2}},{x_3}$. It is important to
 observe that $\set \com x_1,{x_2,x_3}\endset^{T_0}=\set \com x_1,{x_2,x_3}\endset^{T}$.

 We shall also let $S^2$ denote the $T$-space of $\konex$ that is generated by $\com x_1,{x_2}$;
 that is, $S^2=\set \com x_1,{x_2}\endset^S$, and we point out that $\set\com x_1,{x_2}\endset^{S_0}=
 \set \com x_1,{x_2}\endset^S$.

In 1988 (see \cite{Ok}), S. V.  Okhitin proved that for any field $k$ of characteristic zero, the $T$-space 
$CP(M_2(k))$ is finitely based, and he raised the question as to whether $CP(A)$ is finitely based for every 
(unitary) associative algebra $A$ for which $0\ne T(A)\subsetneq CP(A)$.
Then in 2001, V. V. Shchigolev (see \cite{Sh}) showed that for any field of characteristic zero, every $T$-space of $\kzerox$
is finitely based, and it follows from this that every $T$-space of $\konex$ is also finitely based. This more than
answers Okhitin's question (in the affirmative) for fields of characteristic zero.

For a field of characteristic 2, the infinite-dimensional unitary and nonunitary Grassmann algebras are commutative
and thus each has finitely based $T$-space of central polynomials.

We shall show in the following that if $p>2$ and $k$ is an arbitrary field of characteristic $p$, then neither
$CP(\nonunitgrass)$ nor $CP(\unitgrass)$ is finitely based, thus providing a negative answer to
Okhitin's question.

We shall make frequent reference to the following familiar facts.

\begin{lemma}\label{lemma: handy}
 \begin{list}{(\roman{parts})}{\usecounter{parts}}
  \item 
    $\com u,{vw}=\com u,vw+v\com u,w$ for all $u,v,w\in \kzerox$.
  \item 
    $\com u,{vw}=\com u,vw+\com u,wv+\com v,{\com u,w}$ for all $u,v,w\in \kzerox$.
  \item 
     $\com u,{\prod_{i=1}^n v_i}\cong \sum_{i=1}^n 
          \bigl(\mkern3mu \com u,{v_i}\prod_{\substack{j=1\\j\ne i}}^n v_j\mkern3mu\bigr)\mod{T^{(3)}}$ for any positive 
	  integer $n$, and any $u,v_1,v_2,\ldots,v_n\in \kzerox$. 
   \item  
     $\com u,v\com w,x\cong -\com u,w\com v,x\mod{T^{(3)}}$ for all $u,v,w\in \kzerox$.
   \item  
    $\com u,v\com u,w\cong 0\mod{T^{(3)}}$ for all $u,v,w\in \kzerox$.
   \item  
    $\com u,v uw\cong \com u,vwu\mod{T^{(3)}}$ for all $u,v,w\in \kzerox$.    
  \item 
   For any $n\ge2$, $x_1^nx_2^n\cong(x_1x_2)^n+\choice n,2\com x_1,{x_2}x_1^{n-1}x_2^{n-1}\mod{T^{(3)}}$.    
 \end{list}
\end{lemma}

\begin{proof}
 Parts (i) and (ii) are evident by direct calculation. We prove (iii) by induction on $n$. The result is trivial 
 when $n=1$, while  (ii) establishes the case for $n=2$.
 Suppose now that $n\ge2$ is an integer for which the result holds, and consider $u,v_1,v_2,\ldots,v_n,v_{n+1}\in\kzerox$.
 Working modulo $T^{(3)}$, we have
 \begin{align*}
  \com u,{\prod_{i=1}^{n+1} v_i}&\cong\com u,{\prod_{i=1}^n v_i} v_{n+1} + \com u,{v_{n+1}}\prod_{j=1}^n v_i\quad\mbox{by (i)}\\
                                &\cong\sum_{i=1}^n\bigl(\mkern3mu\com u,{v_i}\prod_{\substack{j=1\\j\ne i}}^{n+1}v_i\mkern3mu\bigr)+\com u,{v_{n+1}}\prod_{j=1}^n v_i\quad\text{by the inductive hypothesis}\\
				&=\sum_{i=1}^{n+1}\bigl(\mkern3mu\com u,{v_i}\prod_{\substack{j=1\\j\ne i}}^{n+1}v_i\mkern3mu\bigr).
  \end{align*}

 The result follows now by induction.

 For (iv) and (v), see Latyshev \cite{L}. 
 (vi) follows immediately from (v), as we have $\com u,v uw-\com u,v wu=\com u,v\com u,w\cong 0\mod{T^{(3)}}$.
 
 For (vii), observe that $x_1^2x_2^2=(x_1x_2)^2+x_1\com x_1,{x_2}x_2\cong (x_1x_2)^2+\com x_1,{x_2}x_1x_2$
 modulo $T^{(3)}$. Suppose now that $n\ge 2$ is an integer such that $x_1^nx_2^n\cong(x_1x_2)^n+\choice n,2\com x_1,{x_2}x_1^{n-1}x_2^{n-1}$,
 again modulo $T^{(3)}$.
 Then 
 \begin{align*}
    x_1^{n+1}x_2^{n+1}&=x_1(x_1^nx_2^n)x_2\cong x_1(x_1x_2)^nx_2+\choice n,2x_1\com x_1,{x_2}x_1^{n-1}x_2^{n}\\
     &\cong (x_1x_2)^{n+1}+\com x_1,{(x_1x_2)^n}x_2+\choice n,2\com x_1,{x_2}x_1^nx_2^n\\
     &\cong (x_1x_2)^{n+1}+n\com x_1,{x_1x_2}(x_1x_2)^{n-1}x_2+\choice n,2\com x_1,{x_2}x_1^nx_2^n\quad\text{by (iii)}\\
     &\cong (x_1x_2)^{n+1}+\choice n+1,2\com x_1,{x_2}x_1^nx_2^n\quad\text{by (vi)},
 \end{align*}
 and so the result follows by induction.
\end{proof}

\begin{lemma}\label{lemma: need for alt description of S gen}
 Let $\alpha_1,\alpha_2$ be positive integers, and suppose that $u\in\kzerox$ is central modulo $T^{(3)}$.
 Then $\com x_2,{x_1^{\alpha_1+1}x_2^{\alpha_2}u}\cong (\alpha_1+1)\com x_2,{x_1} x_1^{\alpha_1}x_2^{\alpha_2}u\mod{T^{(3)}}$
 for any $x_1,x_2\in X$.
\end{lemma}

\begin{proof}
 Working modulo $T^{(3)}$, we have
 \begin{align*}
  \com x_2,{x_1^{\alpha_1+1}x_2^{\alpha_2}u}&\cong \com x_2,{x_1^{\alpha_1+1}x_2^{\alpha_2}}u+
              \com x_2,u x_1^{\alpha_1+1}x_2^{\alpha_2} \quad{\text{by Lemma \ref{lemma: handy} (ii)}}\\ 
   &\cong  \com x_2,{x_1^{\alpha_1+1}x_2^{\alpha_2}}u
               \quad{\text{since $\com x_2,u\cong0\mod{T^{(3)}}$}}\\ 
   &\cong  \com x_2,{x_2^{\alpha_2}}x_1^{\alpha_1+1}u+\com x_2,{x_1^{\alpha_1+1}}x_2^{\alpha_2}u \quad{\text{by Lemma \ref{lemma: handy} (ii)}}\\ 
   &\cong (\alpha_1+1)\com x_2,{x_1} x_1^{\alpha_1}x_2^{\alpha_2}u\quad{\text{by Lemma \ref{lemma: handy} (iii).}}
 \end{align*}
\end{proof}


\section{For $p>2$, $CP(\nonunitgrass)$ is not finitely based}\label{nonunitary section}

We assume from this point on that $p>2$.

\begin{definition}\label{definition: siderov's elements}
 Let $\siderovset$ denote the set of all elements of the form 
 \begin{list}{(\roman{parts})}{\usecounter{parts}}
  \item $\prod_{r=1}^t x_{i_r}^{\alpha_r}$, or
  \item $\prod_{r=1}^s \com x_{j_{2r-1}},{x_{2r}}x_{j_{2r-1}}^{\beta_{2r-1}}x_{j_{2r}}^{\beta_{2r}}$, or
  \item $\bigl(\prod_{r=1}^t x_{i_r}^{\alpha_r}\bigr)\prod_{r=1}^s \com x_{j_{2r-1}},{x_{2r}}x_{j_{2r-1}}^{\beta_{2r-1}}x_{j_{2r}}^{\beta_{2r}}$, 
 \end{list}  
 \noindent where $\set i_1,\ldots,i_r\endset\cap\set j_1,\ldots,j_{2s}\endset=\nullset$, $j_1<j_2<\cdots <j_{2s}$, 
 $\beta_i\ge0$ for all $i$, $i_1<i_2<\cdots < i_t$, and $\alpha_i\ge 1$ for all $i$.

 Let $u\in \siderovset$. If $u$ is of the form (i), then the beginning
 of $u$, $\beg(u)$, is $\prod_{r=1}^t x_{i_r}^{\alpha_r}$, the end of
 $u$, $\endof(u)$, is empty, the length of the beginning of $u$,
 $\lbeg{u}$, is equal to $t$ and the length of the end of $u$,
 $\lend(u)$, is 0. If $u$ is of the form (ii), then we say that
 $\beg(u)$, the beginning of $u$, is empty, $\endof(u)$, the end of $u$,
 is $\prod_{r=1}^s \com
 x_{j_{2r-1}},{x_{2r}}x_{j_{2r-1}}^{\beta_{2r-1}}x_{j_{2r}}^{\beta_{2r}}$,
 and $\lbeg{u}=0$ and $\lend(u)=s$. If $u$ is of the form (iii),  then
 we say that $\beg(u)$, the beginning of $u$, is $\prod_{r=1}^t
 x_{i_r}^{\alpha_r}$, $\endof(u)$, the end of $u$, is $\prod_{r=1}^s
 \com
 x_{j_{2r-1}},{x_{2r}}x_{j_{2r-1}}^{\beta_{2r-1}}x_{j_{2r}}^{\beta_{2r}}$,
 and $\lbeg{u}=t$ and $\lend(u)=s$.
 
 Finally, let
  \begin{align}
   \boundedsiderovset&=\{\, u\in \siderovset\mid \text{for each $i$, }\deg_{x_i}(u)<p\text{ if $x_i$ appears in $\beg(u)$} \notag\\
      &\hskip110pt\text{or $\deg_{x_i}(u)\le p$ if $x_i$ appears in $\endof(u)$}\,\}.\notag
  \end{align}
\end{definition}

\begin{lemma}[Theorem 3 of \cite{Si}]\label{lemma: identities for inf dim nonunitary}
 $T(\nonunitgrass)=\set x_1^p\endset^{T_0}+T^{(3)}$.
\end{lemma} 

\begin{definition}\label{definition: def of S}
 For $u,v\in\konex$, let $\kappa(u,v)=\com u,vu^{p-1}v^{p-1}$. Then for each
 $m\ge1$, let $w_m=\prod_{r=1}^m\kappa(x_{2r-1},x_{2r})$. 
\end{definition}

\begin{theorem}[Theorem 1.3 of \cite{Ra}]\label{theorem: main theorem for nonunitary}
  For $k$ any field of characteristic $p>2$, $CP(\nonunitgrass)= S^2+\set w_m\mid \ m\ge 1\endset^S+T(\nonunitgrass)$.
\end{theorem}
  
\begin{lemma}\label{bss li}
 $\set u+T(\nonunitgrass)\mid u\in \boundedsiderovset\endset$ is a linear basis for $\kzerox/T(\nonunitgrass)$.
\end{lemma}

\begin{proof} 
  By Lemma 2.10 of \cite{Si} (or Lemma 2.4 of \cite{Ra}), $\set u+T(\nonunitgrass)\mid u\in \boundedsiderovset\endset$ spans
  $\kzerox/T(\nonunitgrass)$, and Theorem 3 of Siderov \cite{Si} implies that $\gen{\boundedsiderovset}\cap T(\nonunitgrass)=\set 0\endset$.
\end{proof}

\begin{definition}\label{special bss} Let
 $$
  \specialsiderovset=\set u\in\boundedsiderovset\mid \text{ there is $x\in X$ in $\endof(u)$ such that $\deg_{x}(u)<p$}\endset.
 $$ 
\end{definition}

Note that as $\specialsiderovset\subseteq \boundedsiderovset$, it follows from Lemma \ref{bss li} that
$\set u+T(\nonunitgrass)\mid u\in \specialsiderovset\endset$ is linearly independent in $\kzerox/T(\nonunitgrass)$.

\begin{lemma}\label{s2 guys}
 $S^2+T(\nonunitgrass)\subseteq \gen{\specialsiderovset}+T(\nonunitgrass)$.
\end{lemma}

\begin{proof}
 It suffices to prove that for any $u,v\in\kzerox$, $\com u,v\in \gen{\specialsiderovset}+T(\nonunitgrass)$.
 By Lemma 4.1 of \cite{Ch} (we note that while the results of \cite{Ch} were formulated for the case of characteristic zero,
 the proof of Lemma 4.1 is valid in general), for any $u,v\in \kzerox$, $\com u,v$ can be written, modulo $T^{(3)}$, as 
 a sum of terms of the form $\com u_i,{x_i}$, where $u_i$ is a monomial in $\kzerox$ and $x_i\in X$. Now by Lemma 2.10
 of \cite{Si}, modulo $T^{(3)}$, each monomial of $\kzerox$ can be written as a linear combination of elements of 
 $\boundedsiderovset$, so it suffices to prove that modulo $T(\nonunitgrass)$, for each $v\in\boundedsiderovset$ and $x\in X$, $\com v,x$ can be written
 as a linear combination of elements from $\specialsiderovset$. Let $v\in\boundedsiderovset$ and $x\in X$. If $\beg(v)$ is empty,
 Then $\com v,x\in T^{(3)}$ and thus $\com v,x$ is zero modulo $T(\nonunitgrass)$. Suppose now that $\beg(v)$ is not
 empty, so $v=\prod_{r=1}^t x_{i_r}^{\alpha_r}\prod_{i=1}^s\com x_{j_{2i-1}},
 {x_{j_{2i}}}x_{j_{2i-1}}^{\beta_{2i-1}}x_{j_{2i}}^{\beta_{2i}}\in \boundedsiderovset$. By
 Lemma \ref{lemma: handy} (iii) and (vi), $\com v,x$ can be written as a linear combination
 of terms each of the form $d_j=\bigl(\prod_{\substack{r=1\\r\ne j}}^t x_{i_r}^{\alpha_r}\bigr)x_{i_j}^{\alpha_j-1}\com x_{i_j},x
 \prod_{i=1}^s\com x_{j_{2i-1}},{x_{j_{2i}}}x_{j_{2i-1}}^{\beta_{2i-1}}x_{j_{2i}}^{\beta_{2i}}$. 
 By Lemma \ref{lemma: handy} (iv) and (v), for each $j$, modulo $T(\nonunitgrass)$,
 $d_j$ is congruent either to 0 or (up to sign) to an element $l_j$ of $\specialsiderovset$ since $x_{i_j}$ has degree at most $p-1$ in $d_j$
 and appears in $\endof(d_j)$, hence in $\endof(l_j)$.
\end{proof} 

\begin{corollary}\label{corollary: fundamental fact}
 For each $m\ge1$, $w_{m}\notin S^2+T(\nonunitgrass)$.
\end{corollary}

\begin{proof}
 Let $m\ge1$. Then $w_m\in \boundedsiderovset$,
 and $w_m\notin \specialsiderovset$, so by Lemma \ref{bss li}, $\set w_m\endset\cup \specialsiderovset$ is linearly
 independent modulo $T(\nonunitgrass)$. Thus $w_m\notin \gen{\specialsiderovset}+T(\nonunitgrass)$, hence
 by Lemma \ref{s2 guys}, $w_{m}\notin S^2+T(\nonunitgrass)$.
\end{proof}
 
In applications of the following result, it will be important to recall that $S^2$ is the $T$-space generated
by $\com x_1,{x_2}$ in either $\konex$ or $\kzerox$, and $T^{(3)}$ is the $T$-ideal generated by $\com x_1,{x_2,x_3}$
in either $\konex$ or $\kzerox$.

\begin{lemma}\label{w additive}
 For any $u,v,w\in\konex$, modulo $T^{(3)}$, we have
 $$
  \kappa(u,v\mkern-3mu +\mkern-3mu w)\cong \kappa(u,v)\mkern-2mu + \mkern-2mu\kappa(u,w)\mkern-2mu +\mkern-2mu
 \sum_{i=0}^{p-2}(i\mkern-3mu + \mkern-3mu1)^{-1} \choice p-1,i \com u,{v^{i+1}w^{p-(i+1)}u^{p-1}}.
 $$ 
\end{lemma} 

\begin{proof}
 Recall that $\kappa(u,v+w)=\com u,{v+w}u^{p-1}(v+w)^{p-1}$. To begin with, we shall prove that 
 \begin{align*}
 \com u,{v+w}(v+w)^{p-1}& \cong \com u,v v^{p-1}+\com u,w w^{p-1}\\
 &\hskip25pt+\sum_{i=0}^{p-2}(i\mkern-3mu + \mkern-3mu1)^{-1} \choice p-1,i \com u,{v^{i+1}w^{p-(i+1)}}\mod{T^{(3)}}.
 \end{align*}
  
 \noindent We have
 \begin{align*}
  \com u,{v+w}(v+w)^{p-1}&=\com u,v (v+w)^{p-1}+\com u,w (v+w)^{p-1}\\
   &\mkern-6mu\overset{T^{(3)}}{\cong} \com u,v v^{p-1}+\com u,v\sum_{i=0}^{p-2}\choice p-1,i v^iw^{p-1-i}\\
   &\hskip40pt+\com u,w w^{p-1} +\com u,w\sum_{i=1}^{p-1}\choice p-1,i v^iw^{p-1-i}\\
 \end{align*}  
 
 \noindent so it suffices to show that $\com u,v\sum_{i=0}^{p-2}\choice p-1,i v^iw^{p-1-i}
   +\com u,w\sum_{i=1}^{p-1}\choice p-1,i v^iw^{p-1-i}$ is congruent to
 $\sum_{i=0}^{p-2}(i\mkern-3mu + \mkern-3mu1)^{-1} \choice p-1,i \com u,{v^{i+1}w^{p-(i+1)}}$ modulo $T^{(3)}$.
 By Lemma \ref{lemma: need for alt description of S gen}, we have 
 $$
    \com u,v\sum_{i=0}^{p-2}\choice p-1,i v^iw^{p-1-i}\cong
    \sum_{i=0}^{p-2}(i+1)^{-1}\com u,{v^{i+1}}\choice p-1,i w^{p-1-i}\mod{T^{(3)}},
 $$
 and by Lemma \ref{lemma: handy} (ii), $\sum_{i=0}^{p-2}(i+1)^{-1}\com u,{v^{i+1}}\choice p-1,i w^{p-1-i}$
 is congruent modulo $T^{(3)}$ to
 $$
  \sum_{i=0}^{p-2}(i+1)^{-1}\com u,{v^{i+1}w^{p-1-i}}\choice p-1,i\\
     -\sum_{i=0}^{p-2}(i+1)^{-1}\com u,{w^{p-1-i}}\choice p-1,iv^{i+1}.
 $$
 It is sufficient therefore to prove that 
 $$
  \sum_{i=0}^{p-2}(i+1)^{-1}\com u,{w^{p-1-i}}\choice p-1,iv^{i+1}\cong
    \com u,w\sum_{i=1}^{p-1}\choice p-1,i v^iw^{p-1-i}\mod{T^{(3)}}.
 $$
 But by Lemma \ref{lemma: need for alt description of S gen}, together with a change of variable, we have 
 $$
  \com u,w\sum_{i=1}^{p-1}\choice p-1,i v^iw^{p-1-i}
    \cong \sum_{i=0}^{p-2}\choice p-1,{i+1} (p-i-1)^{-1} \com u,{w^{p-i-1}}v^{i+1}.
 $$
 Since $p-i-1=-(i+1)$ and for each $i$ with $0\le i\le p-2$, $0=\choice p,{i+1}=\choice p-1,i+\choice p-1,{i+1}$,
 and thus $\choice p-1,{i+1}=-\choice p-1,i$, the result follows. 
    
 To complete the proof, observe that by Lemma \ref{lemma: handy} (ii), for each $i$,
 $$
    \com u,{v^{i+1}w^{p-(i+1)}}u^{p-1}\cong \com u,{v^{i+1}w^{p-(i+1)}u^{p-1}}\mod{T^{(3)}}.
 $$
 \end{proof}
 
 The following additivity result is fundamental for this work.
 
 \begin{corollary}\label{w general additive}
  For any $m\ge1$, let $u_1,u_2,\ldots,u_{2m},v\in\konex$. Then modulo $S^2+T^{(3)}$, for any $i$ with
  $1\le i\le 2m$, 
  \begin{align*}
    w_m(u_1,u_2,\ldots,u_i+v,\ldots,u_{2m})&\cong w_m(u_1,u_2,\ldots,u_i,\ldots,u_{2m})\\
    &\hskip30pt+w_m(u_1,u_2,\ldots,v,\ldots,u_{2m}).
  \end{align*}
 \end{corollary}
 
 \begin{proof}
  First we note that for any $u,v\in\konex$, $\kappa(u,v)$ is central modulo $T^{(3)}$. Moreover, $\kappa(v,u)\cong
  -\kappa(u,v)\mod{T^{(3)}}$, so it suffices to prove the result for even $i$. Thus we shall consider $1\le i\le m$,
  and let $\gamma=\prod_{\substack{j=1\\j\ne i}}^m\kappa(u_{2j-1},u_{2j})$. Then $\gamma$ is central modulo $T^{(3)}$,
  and so 
  $$
   w_m(u_1,u_2,\ldots,u_{2i}+v,\ldots,u_{2m})\cong \kappa(u_{2i-1},u_{2i}+v)\gamma\mod{T^{(3)}}.
  $$
  By Lemma \ref{w additive}, 
  \begin{align*}
   \kappa(u_{2i-1},u_{2i}+v) &\cong \kappa(u_{2i-1},u_{2i}) +\kappa(u_{2i-1},v)\\
  &\hskip3pt+\sum_{r=0}^{p-2} (r\mkern-3mu + \mkern-3mu1)^{-1}\choice p-1,r \com u_{2i-1},{u_{2i}^{r+1}v^{p-(r+1)}u_{2i-1}^{p-1}}\mod{T^{(3)}},
  \end{align*}
  and thus
  $w_m(u_1,u_2,\ldots,u_{2i}+v,\ldots,u_{2m})\cong \kappa(u_{2i-1},u_{2i})\gamma +\kappa(u_{2i-1},v)\gamma
  +\sum_{r=0}^{p-2}(r\mkern-3mu + \mkern-3mu1)^{-1}\choice p-1,r \com u_{2i-1},{u_{2i}^{r+1}v^{p-(r+1)}u_{2i-1}^{p-1}}\gamma\mod{T^{(3)}}$.
  Finally, by Lemma \ref{lemma: handy} (ii) and the fact that $\gamma$ is central modulo $T^{(3)}$, we have for
  each $r$ that $\com u_{2i-1},{u_{2i}^{r+1}v^{p-(r+1)}u_{2i-1}^{p-1}}\gamma\cong
  \com u_{2i-1},{u_{2i}^{r+1}v^{p-(r+1)}u_{2i-1}^{p-1}\gamma}\mod{T^{(3)}}$. It follows now 
  that
  $\sum_{r=0}^{p-2}(r\mkern-3mu + \mkern-3mu1)^{-1}\choice p-1,r \com u_{2i-1},{u_{2i}^{r+1}v^{p-(r+1)}u_{2i-1}^{p-1}}\gamma\in S^2+T^{(3)}$, as required.
 \end{proof}
 
 Thus for any $m\ge1$, modulo $S^2+T^{(3)}$, $w_m$ is additive in each variable $x_1,x_2,\ldots,x_{2m}$. 
 
 \begin{lemma}\label{multiplicative property of w}
   Let $u,v,w\in\konex$. Then 
   $$
    \kappa(u,vw)\cong v^p\kappa(u,w)+w^p\kappa(u,v)\mod{T^{(3)}}.
   $$ 
   In particular, if $u,v\in\kzerox$ and $\alpha\in k$, then $\kappa(u,\alpha v)=\alpha^p\kappa(u,v)$.
 \end{lemma} 
 
 \begin{proof}
  By Lemma \ref{lemma: handy} (ii), $\com u,vw\cong \com u,vw+\com u,wv\mod{T^{(3)}}$, so we have
  \begin{align*}
  \kappa(u,vw)&=\com u,{vw}u^{p-1}(vw)^{p-1}\\
  &\cong \com u,vwu^{p-1}(vw)^{p-1}+\com u,wvu^{p-1}(vw)^{p-1}\mod{T^{(3)}}
  \end{align*}
  By Lemma \ref{lemma: handy} (vi), in any product expression with $\com u,v$ and $u$ as factors, $u$ commutes within the
  product expression, modulo $T^{(3)}$.
  Thus
  \begin{align*}
   \com u,vwu^{p-1}(vw)^{p-1}&\cong \com u,vwu^{p-1}v^{p-1}w^{p-1}\cong \com u,vu^{p-1}v^{p-1}w^{p}\\
   &=\kappa(u,v)w^p\mod{T^{(3)}}.
  \end{align*} 
  Since $\kappa(u,v)$ is central modulo $T^{(3)}$, we have 
  $$
   \com u,vwu^{p-1}(vw)^{p-1}\cong w^p\kappa(u,v)\mod{T^{(3)}}.
  $$
  Similarly, $\com u,wvu^{p-1}(vw)^{p-1}\cong v^p\kappa(u,w)\mod{T^{(3)}}$.
 \end{proof} 
 
 \begin{corollary}\label{w multiplicative}
  For any $u,v,w\in\kzerox$, $\kappa(u,(vw))\cong 0\mod{T(\nonunitgrass)}$.
 \end{corollary}
 
 \begin{proof}
   This follows immediately from Lemma \ref{multiplicative property of w} since $x_1^p\in T(\nonunitgrass)$.
 \end{proof}
 
 Let $\zplus$ denote the set of positive integers.
 
 \begin{definition}\label{w sets}
  For each $m\ge1$, let $I_m$ denote the set of all strictly increasing functions from $J_{2m}=\set 1,2,\ldots,2m\endset$
  into $\zplus$ (that is, $f(i)<f(j)$ if $i<j$), and let $W_m=\set w_j(x_{f(1)},x_{f(2)},\ldots,x_{f(2j)})
  \mid 1\le j\le m,\ f\in I_j\endset$. Finally, let $W=\bigcup_{m=1}^\infty W_m$.
 \end{definition}

 \begin{lemma}\label{chain of t spaces}
  Suppose that $V_i$, $i\ge1$ are $T$-spaces of $\kzerox$ (or $\konex$) such that for each $i$, $V_i\subsetneq V_{i+1}$. Then
  $V=\bigcup_{i=1}^\infty V_i$ is a $T$-space of $\kzerox$ (respectively, $\konex$) that is not finitely based.
 \end{lemma}
 
 \begin{proof}
  If $V$ were finitely based, then for some $n$, $V_n$ would contain a basis for $V$, and thus $V=V_n\subsetneq V_{n+1}\subseteq V$,
  which is not possible.
\end{proof}

 \begin{lemma}\label{basic fact}
  $\set u+S^2+T(\nonunitgrass)\mid u\in W\endset$ is a linear basis for the vector space $CP(\nonunitgrass)/(S^2+T(\nonunitgrass))$.
 \end{lemma}
 
 \begin{proof}
  Since $CP(\nonunitgrass)=W^S+S^2+T(\nonunitgrass))$, by Corollary \ref{w general additive}, it suffices to consider only 
  linear combinations of elements of the form
  $w_m(u_1,u_2,\ldots,u_{2m})$, $u_i\in \kzerox$, and by Corollary \ref{w multiplicative}, it suffices to consider only
  such elements where $u_i\in X$ for each $i$. For any subset of size $2m$ in $\zplus$, say
  $\set i_1,i_2,\ldots,i_{2m}\endset$, there exists $\sigma\in S_{2m}$ such that 
  $\sigma(i_1)<\sigma(i_2)<\cdots <\sigma(i_{2m})$, and
  by Lemma \ref{lemma: handy} (v), $w_m(x_{i_1},\ldots,x_{i_{2m}})\cong(-1)^{\text{sgn}(\sigma)}w_m(x_{\sigma(i_1)},\ldots,x_{\sigma(i_{2m})})\mod{T^{(3)}}$.
  This proves that $\set u+S^2+T(\nonunitgrass)\mid u\in W\endset$ spans $CP(\nonunitgrass)/(S^2+T(\nonunitgrass))$.
  
  In order to establish linear independence, suppose that $u\in S^2+T(\nonunitgrass)$ is a linear combination of elements of $W$.
  Order the set $I=\bigcup_{m=1}^\infty I_m$ lexically (so that for $m_1<m_2$, $f_1\in I_{m_1}$, and $f_2\in I_{m_2}$, we have
  $f_1<f_2$). Then there exists a smallest $f\in I$ such that for some nonzero $\alpha\in k$, $\alpha w_m(x_{f(1)},\ldots,x_{f(2m)})$
  is a summand of $u$.
  Let $\theta$ denote the endomorphism of $\kzerox$ that is determined by mapping $x_{f(i)}$ to $x_i$
  for each $i=1,2,\ldots,2m$, and mapping all other elements of $X$ to 0. Since $S^2+T(\nonunitgrass)$ is a $T$-space,
  $\alpha w_m=\theta(u)\in S^2+T(\nonunitgrass)$. But by Corollary \ref{corollary: fundamental fact}, $w_m\notin S^2+T(\nonunitgrass)$,
  so $\alpha=0$. Since $\alpha \ne0$, we have a contradiction and thus the linear independence is established.
 \end{proof}

%
 
\begin{theorem}\label{cp nonunitary not fb}
 For any prime $p>2$, and any field of characteristic $p$, the $T$-space $CP(\nonunitgrass)$ is not finitely based.
\end{theorem}

\begin{proof}
 We have $CP(\nonunitgrass)=W^S+S^2+T(\nonunitgrass)$. For each $m\ge1$, let $U_m=W_m^S+S^2+T(\nonunitgrass)$.
 Then $CP(\nonunitgrass)=\bigcup_{m=1}^\infty U_m$, and for each $m\ge1$, $U_m\subsetneq U_{m+1}$, where the
 inequality follows from Lemma \ref{basic fact}. It follows now from Lemma \ref{chain of t spaces} that
 $CP(\nonunitgrass)$ is not finitely based.
\end{proof} 
   
 The following result is the nonunitary analogue of \cite{Shch}, Theorem 4.

\begin{corollary}\label{shchigolev analogue}
 The $T$-space $W^S+T(\nonunitgrass)$ is not finitely based.
\end{corollary}

\begin{proof}
 If $W^S+T(\nonunitgrass)$ is finitely based, then $CP(\nonunitgrass)=W^S+T(\nonunitgrass)+S^2$ is finitely
 based, which contradicts Theorem \ref{cp nonunitary not fb}.
\end{proof}

\begin{corollary}
 $W^S$ is not finitely based.
\end{corollary}

\begin{proof}
 Since $T(\nonunitgrass)$ is finitely based, this follows immediately from Corollary \ref{shchigolev analogue}.
\end{proof}
  

\section{For $p>2$, $CP(\unitgrass)$ is not finitely based}
 We extend the definition of $w_m$ by setting $w_0=1$.
 It was shown in \cite{Ra} that $CP(\unitgrass)=S^2+\set x_0^pw_m\mid m\ge0\endset^S+T(\unitgrass)$ if $k$ is an
 infinite field of characteristic $p>2$. Subsequently, we showed in \cite{ChuRa} that the same is true even if the field is
 finite. The difference between the two situations is in the expression for $T(\unitgrass)$. If $k$ is infinite,
 then it was shown in \cite{Gi} that $T(\unitgrass)=T^{(3)}$, while if $k$ is finite, say of size $q$, then it was shown
 in \cite{ChuRa} that $T(\unitgrass)=\set x_1^{qp}-x_1^p\endset^T+T^{(3)}$.
 
 \begin{lemma}\label{w additive in unitary}
  Let $m\ge1$, and let $\alpha_1,\alpha_2,\ldots,\alpha_{2m}\in k$. Then
  $$
   w_m(x_1+\alpha_1,x_2+\alpha_2,\ldots,x_{2m}+\alpha_{2m})\cong w_m\mod{S^2+T^{(3)}}.
   $$
 \end{lemma}
 
 \begin{proof}
  By Corollary \ref{w general additive}, modulo $S^2+T^{(3)}$, $w_m(u_1,u_2,\ldots,u_i+v,\ldots,u_{2m})$
  is congruent to
  $$
   w_m(u_1,u_2,\ldots,u_i,\ldots,u_{2m})  +w_m(u_1,u_2,\ldots,v,\ldots,u_{2m})
  $$
  for any $u_1,u_2,\ldots,u_{2m},v\in\konex$. Since $\kappa(u,v)=0$ if $u$ or $v$ is an element of $k$,
  it follows that modulo $S^2+T^{(3)}$, we have
  $$
   w_m(x_1+\alpha_1,x_2+\alpha_2,\ldots,x_{2m}+\alpha_{2m})\cong w_m(x_1,x_2,\ldots,x_{2m})=w_m
  $$
 \end{proof}
 
 \begin{lemma}\label{rep of cp}
  $CP(\unitgrass)=S^2+\set x_0^pw_{m}\mid m\ge0\endset^{S_0}+\set w_{m}\mid m\ge 0\endset^S+T(\unitgrass)$.
 \end{lemma}
 
 \begin{proof}
  Let $U=S^2+\set x_0^pw_{m}\mid m\ge0\endset^{S_0}+\set w_{m}\mid m\ge 0\endset^S+T(\unitgrass)$.
  Since $CP(\unitgrass)=S^2+\set x_0^pw_{m}\mid m\ge0\endset^S+T(\unitgrass)$, it is evident
  that $U\subseteq CP(\unitgrass)$. It remains to prove that 
  $\set x_0^pw_{m}\mid m\ge0\endset^S\subseteq U$.
  Let $u\in \kzerox$ and $\alpha\in k$. Then $(u+\alpha)^p=u^p+\alpha^p$, $\alpha^p\in k\subseteq\set w_{m}\mid m\ge0\endset^S$,
  and $u^p\in\set x_0^pw_{m}\mid m\ge0\endset^{S_0}$, so $(u+\alpha)^p\in U$. Next, let $m\ge 1$ and let $u,u_1,u_2,\ldots,u_{2m}\in\kzerox$ and
  $\alpha_1,\ldots,\alpha_{2m}\in k$. By Lemma \ref{w additive in unitary}, there is $v\in S^2$ such that
  $w_m(u_1+\alpha_1,\ldots,u_{2m}+\alpha_{2m})  \cong w_m(u_1,\ldots,u_{2m})+v\mod{T(\unitgrass)}$. As
  $(u+\alpha)^p=u^p+\alpha^p$, we have $(u+\alpha)^pw_m(u_1+\alpha_1,\ldots,u_{2m}+\alpha_{2m})\cong
  u^pw_m(u_1,\ldots,u_{2m})+\alpha^pw_m(u_1,\ldots,u_{2m})+(u+\alpha)^pv\mod{T(\unitgrass)}$. Now, since
  $x^p$ is a central polynomial for $\unitgrass$, $(u+\alpha)^pv\in S^2+T(\unitgrass)$ (by Lemma \ref{lemma: handy} (ii),
  for any $u,v,w\in\konex$, $\com u,{vw}\cong \com u,vw+\com u,wv
  \mod{T^{(3)}}$, and if $v$ is a central polynomial of $\unitgrass$,
  then $\com u,v\in T(\unitgrass)$). Thus $(u+\alpha)^pw_m(u_1+\alpha_1,\ldots,u_{2m}+\alpha_{2m})\cong
  u^pw_m(u_1,\ldots,u_{2m})+\alpha^pw_m(u_1,\ldots,u_{2m})\mod{S^2+T(\unitgrass)}$, and so 
  $$
   (u+\alpha)^pw_m(u_1+\alpha_1,\ldots,u_{2m}+\alpha_{2m})\in U.
  $$ 
 \end{proof}  
 
 \begin{lemma}\label{wm not in s2 +xp +t3}
  For every $m\ge0$, $w_m\notin S^2+\set x_0^pw_{j}\mid j\ge0\endset^{S_0}+T(\unitgrass)$.
 \end{lemma}
 
 \begin{proof}
  First, note that $S^2+\set x_0^pw_{j}\mid j\ge0\endset^{S_0}+T(\unitgrass)\subseteq S^2+\set x_0^p\endset^{T_0}+T(\unitgrass)$.
  Now, $T(\unitgrass)=T^{(3)}$ if $k$ is infinite, while $T(\unitgrass)=\set(x^{qp}-x^p\endset^T+T^{(3)}$ if $k$ is
  finite of size $q$, and in either case, $T(\nonunitgrass)=\set x_o^p\endset^{T_0}+T^{(3)}$. Thus if $k$ is infinite, 
  we have $S^2+\set x_0^pw_{j}\mid j\ge0\endset^{S_0}+T(\unitgrass)\subseteq S^2+T(\nonunitgrass)$.
  Suppose now that $k$ is finite. As shown in Section 2 of \cite{CR}, for any $\alpha\in k$ and $u\in\kzerox$, we have
  $(\alpha+u)^{qp}-(\alpha+u)^p=u^{qp}-u^p$, so $\set x_0^{qp}-x_0^p\endset^T=\set x_0^{qp}-x_0^p\endset^{T_0}\subseteq \set x_0^p\endset^{T_0}$.
  Thus $S^2+\set x_0^pw_{j}\mid j\ge0\endset^{S_0}+T(\unitgrass)\subseteq S^2+T(\nonunitgrass)$ in this case as well.
  
  The result follows now from Corollary \ref{corollary: fundamental fact}.
 \end{proof}
 
 \begin{lemma}\label{wm multiplicative}
  Let $m\ge1$. Then for any $i$ with $1\le i\le 2m$, 
  $$
   w_m(x_1,x_2,\ldots,x_ix_{2m+1},\ldots,x_{2m})\in \set x_0^pw_{j}\mid j\ge0\endset^{S_0}+T(\unitgrass),
  $$ 
  while for any $\alpha\in k$,
  $$
   w_m(x_1,x_2,\ldots,\alpha x_i,\ldots,x_{2m})=\alpha^pw_m(x_1,x_2,\ldots,x_i,\ldots,x_{2m}).
  $$ 
 \end{lemma}
 
 \begin{proof}
  By Lemma \ref{lemma: handy} (vi), $w_1(x_1,x_2)\cong -w_1(x_2,x_1)\mod{T^{(3)}}$, and for any $m\ge1$,
  $w_m\in CP(\unitgrass)$, so without loss of generality, it suffices to prove the result for $i=2m$. 
  If $m=1$, then we have by Lemma \ref{multiplicative property of w} that 
  $w_1(x_1,x_2x_3)\cong x_2^pw_1(x_1,x_3)+x_3^pw_1(x_1,x_2)\mod{T^{(3)}}$
  and so the result holds in this case. Suppose now that $m>1$. Since 
  $$
   w_m(x_1,x_2,\ldots,x_{2m}x_{2m+1})=w_{m-1}w_1(x_{2m-1},x_{2m}x_{2m+1}),
  $$
  and $w_{m-1}w_1(x_{2m-1},x_{2m}x_{2m+1})$ is congruent to
  $$
    w_{m-1}x_{2m}^pw_1(x_{2m-1},x_{2m+1})+w_{m-1}x_{2m+1}^pw_1(x_{2m-1},x_{2m})\mod{T^{(3)}},
  $$
  the first assertion follows. The second assertion is obvious.
 \end{proof}  

  Recall that $W$ was introduced in Definition \ref{w sets}.

 \begin{lemma}\label{wm lin ind in unitary case}
  The vector space $CP(\unitgrass)/(S^2+\set x_0^pw_{m}\mid m\ge0\endset^{S_0}+T(\unitgrass))$ has linear basis
  $\set u+S^2+\set x_0^pw_{m}\mid m\ge0\endset^{S_0}+T(\unitgrass)\mid u\in \set 1\endset\cup W\endset$.
 \end{lemma}
 
 \begin{proof}
  By Lemma \ref{rep of cp}, $CP(\unitgrass)$ is equal to 
  $$
       k+W^S +S^2+\set x_0^pw_{m}\mid m\ge0\endset^{S_0}+T(\unitgrass).
  $$ 
  Let $U=S^2+\set x_0^pw_{m}\mid m\ge0\endset^{S_0}+T(\unitgrass)\subseteq \kzerox$.
  The spanning argument in the proof of Lemma \ref{basic fact} is applicable here, with
  the respective roles of Corollary \ref{w general additive} and Corollary \ref{w multiplicative} being played
  by Lemma \ref{w additive in unitary} and Lemma \ref{wm multiplicative}.
  
  Now for linear independence, suppose that $u\in U$ is a linear combination of elements of $\set 1\endset\cup W$.
  Since $U\subseteq \kzerox$, $u$ must be a linear combination of elements of $W$. Then, just as in the proof of Lemma \ref{basic fact}, we order the set 
  $I=\bigcup_{m=1}^\infty I_m$ lexically, and find the smallest $f\in I$ such that 
  for some nonzero $\beta\in k$, $\beta w_m(x_{f(1)},\ldots,x_{f(i_{2m})})$, 
  is a summand of $u$. Let $\theta$ denote the endomorphism of $\konex$ that is determined by mapping 
  $x_{f(i)}$ to $x_i$ for each $i=1,2,\ldots,2m$, and mapping all other elements of $X$ to 0. Since 
  $U$ is a $T$-space, $\beta w_m=\theta(u)\in U$. But by Lemma \ref{wm not in s2 +xp +t3}, $w_m\notin U$,
  so $\beta=0$.  Since $\beta \ne0$, we have a contradiction and thus the linear independence is established.
 \end{proof}  


We are thus able to obtain the unitary analogues of the main results of Section \ref{nonunitary section}.
Let $W'_0=\set 1\endset$, and for every $m\ge1$, let $W'_m=W'_0\cup W_m$. Finally, let $W'=\bigcup_{j=0}^\infty W'_m$.

%

\begin{theorem}\label{cp unitary not fb}
 For any prime $p>2$, and any field of characteristic $p$, the $T$-space $CP(\unitgrass)$ is not finitely based.
\end{theorem}

\begin{proof}
 By Lemma \ref{rep of cp}, $CP(\unitgrass)$ is equal to 
 $$
   (W')^S +S^2+\set x_0^pw_{m}\mid m\ge0\endset^{S_0}+T(\unitgrass).
 $$ 

 For each $m$, let $U_m=(W'_m)^S+S^2+\set x_0^pw_{m}\mid m\ge0\endset^{S_0}+T(\unitgrass)$.
 Then $CP(\unitgrass)=\bigcup_{m=1}^\infty U_m$, and for each $m\ge1$, $U_m\subsetneq U_{m+1}$, where the
 inequality follows from Lemma \ref{wm lin ind in unitary case}. It follows now from Lemma \ref{chain of t spaces} that
 $CP(\unitgrass)$ is not finitely based.
\end{proof} 
   
 The following result is basically Theorem 4 of \cite{Shch}, extended in the sense that
 it holds for all fields of characteristic $p>2$, not just infinite fields.

\begin{corollary}[Theorem 4, \cite{Shch}]\label{shchigolev theorem}
 The $T$-space $(W')^S+T(\unitgrass)$ is not finitely based.
\end{corollary}

\begin{proof}
 If $(W')^S+T(\unitgrass)$ is finitely based, say with finite basis $B$, then $B\cup x_0^pB\cup\set \com x_1,{x_2}\endset$
 is a finite basis for $(W')^S+T(\unitgrass)+\set x_0^pw_{m}\mid m\ge0\endset^{S_0} +S^2=CP(\unitgrass)$, which contradicts
 Theorem \ref{cp unitary not fb}.
\end{proof}

\begin{corollary}
 $(W')^S$ is not finitely based.
\end{corollary}

\begin{proof}
 Since $T(\unitgrass)$ is finitely based, this follows immediately from Corollary \ref{shchigolev theorem}. 
\end{proof}

In \cite{Shch}, Shchigolev introduces polynomials $\varphi'_m=\prod_{j=1}^m x_{2i-1}^{p-1}x_{2i}x_{2i-1}x_{2i}^{p-1}$,
 $m\ge1$, and he proves (essentially) the following result:
 
 \begin{corollary}[Theorem 5, \cite{Shch}]
  The $T$-space of $\konex$ that is generated by the set $\set \varphi'_m\mid m\ge1\endset$ is not finitely based.
\end{corollary}  

\begin{proof}
Let $m\ge1$ be given, let $J_{2m}=\set 1,2,\ldots,2m\endset$ and let $J_{2m}^{\text{odd}}$ denote the odd elements of
$J_{2m}$. For any $I\subseteq J_{2m}^{odd}$, $I\ne \nullset$, we have $I=\set i_1,i_2,\ldots,i_k\endset$ with
$i_1<i_2<\cdots<i_k$, and we shall let $w_{|I|}(x_i,x_{i+1}\mid i\in I)$ denote 
$$
 w_{|I|}(x_{i_1},x_{i_1+1},x_{i_2},x_{i_2+1},\ldots,x_{i_k},x_{i_k+1}).
$$ 
It is easily seen now that for each $m\ge1$, $\varphi_m'$ is congruent modulo $T(\unitgrass)$ to
$$
 w_m+\prod_{i=1}^m x_{2i-1}^px_{2i}^p +
\sum_{\substack{I\subseteq J_{2m}^{\text{odd}}\\ I\ne \nullset,J_{2m}^{\text{odd}}}}
 \bigl(\prod_{i\mkern2mu\in J_{2m}^{\text{odd}}-I}x_i^px_{i+1}^p\bigr) w_{|I|}(x_i,x_{i+1}\mid i\in I).
 $$
For each $m\ge1$, let $R_m=\set \varphi_i'\mid 1\le i\le m\endset^S$ and $R=\cup_{i=1}^\infty R_m$. It follows now by
a straightforward inductive argument, using
the representation of $\varphi_m'$ above, that for every $m\ge1$, $(W_m')^S+(x_0^pW'_m)^{S}+T(\unitgrass)=
R_m+(x_0^pR_m)^{S}+\set x_0^p\endset^S+T(\unitgrass)$, and thus $(W')^S+(x_0^pW')^{S}+T(\unitgrass)=
R+(x_0^pR)^{S}+\set x_0^p\endset^{S}+T(\unitgrass)$ (one shows that $\varphi_m'-w_m$ belongs to both 
$(W_m')^S+(x_0^pW'_m)^{S}+T(\unitgrass)$ and $R+(x_0^pR)^{S}+\set x_0^p\endset^{S}+T(\unitgrass)$
by using the inductive hypothesis to show that $\sum_{\substack{I\subseteq J_{2m}^{\text{odd}}\\ I\ne \nullset,J_{2m}^{\text{odd}}}}\bigl(\prod_{i\mkern2mu\in J_{2m}^{\text{odd}}-I}x_i^px_{i+1}^p\bigr) w_{|I|}(x_i,x_{i+1}\mid i\in I)$
belongs to $(W_m')^S+(x_0^pW'_m)^{S}+T(\unitgrass)$ and to $R+(x_0^pR)^{S}+\set x_0^p\endset^{S}+T(\unitgrass)$, while
of course, $x_1^p\cdots x_{2m}^p$ belongs to both sets).
If $R$ is finitely based, then $R+(x_0^pR)^{S}+\set x_0^p\endset^{S}+T(\unitgrass)$
and thus $(W')^S+(x_0^pW')^{S}+T(\unitgrass)$ is finitely based as well. But then 
\begin{align*}
  CP(\unitgrass)&= (W')^S +S^2+\set x_0^pw_{m}\mid m\ge0\endset^{S_0}+T(\unitgrass)\\
  &= (W')^S +S^2+\set x_0^pw_{m}\mid m\ge0\endset^{S}+T(\unitgrass)\quad\text{by Lemma \ref{w additive in unitary}}\\  
  &=(W')^S+S^2+(x_0^pW')^{S}+T(\unitgrass)
\end{align*}
is finitely based, contradicting Theorem \ref{cp unitary not fb}. Thus $R$ is not finitely based.  
\end{proof}


\begin{thebibliography}{00}
\bibitem{Ra} Chuluundorj Bekh-Ochir and S. A. Rankin, {\em The central polynomials of the 
infinite dimensional unitary and nonunitary Grassmann algebras}, preprint, arXiv:0905.1049

\bibitem{ChuRa} Chuluundorj Bekh-Ochir and S. A. Rankin, {\em The identities and the central polynomials of the
 infinite dimensional unitary Grassmann algebras over a finite field}, preprint, arXiv:0905.1049.

\bibitem{CR} Chuluundorj Bekh-Ochir and S. A. Rankin, {\em The identities and central polynomials of the finite dimensional unitary 
  Grassmann algebras over a finite field}, preprint, arXiv:0905.2388.

\bibitem{Ch} Chuluundorj Bekh-Ochir and D. M. Riley, {\em On the Grassmann $T$-space}, Journal of Algebra and its Applications, 
      {\bf 7} (2008), no. 3, 319--336.

\bibitem{Gi} A. Giambruno and P. Koshlukov, {\em On the identities of the Grassmann algebras in characteristic $p>0$}, Isr. J. Math.
{\bf 122} (2001), 305--316.

\bibitem{L} V. N. Latyshev, {\em On the choice of a basis in a $T$-ideal}, Siberian  Math.  J. {\bf 4} (1963), no. 5, 1122--1127.

\bibitem{Ok} S. V. Okhitin, {\em Central polynomials of an algebra of second-order matrices} (Russian), Vestnik Moskov. Univ. Ser. I Mat. Mekh., (1988), no. 4, 
    61--63; translation in Moscow Univ. Math. Bull., {\bf 43} (1988), no. 4, 49--51. 

\bibitem{Si} Plamen Zh. Chiripov and Plamen N. Siderov, 
   {\em On bases for identities of some varieties of associative algebras} 
   (Russian), PLISKA Stud. Math. Bulgar, {\bf 2} (1981), 103-115.

\bibitem{Sh} V. V. Shchigolev, {\em Finite basis property of {$T$}-spaces over fields of characteristic zero}, Izv. Ross. Akad. Nauk Ser. Mat. {\bf 65}, (2001), no. 5, 191-224; translation in Izv. Math., {\bf 65} (2001), 
   no. 5, 1041-1071.

\bibitem{Shch} V. V. Shchigolev, {\em Examples of $T$-spaces with an infinite basis}, Sbornik Mathematics {\bf 191} (2000), no. 3. 459--476.

\end{thebibliography}
\end{document}